\documentclass[12pt]{amsart}
\usepackage{amsmath,amsthm, amssymb}
\usepackage{subcaption}
\usepackage{graphicx} \graphicspath{{figs/}}
\usepackage{hyperref}
\usepackage{xcolor}
\usepackage[percent]{overpic}
\usepackage{rotating}

\def\MS{M}

\def\R{\mathbb{R}}

\def\Z{\mathbb{Z}} 
 
\def\HH{\mathbb{H}} 
 
\def\Q{\mathbb{Q}} 
\def\QQ{\mathbb{Q}} 
\def\RR{\mathbb{R}}
 
\def\ZZ{\mathbb{Z}} 
\newcommand\bs{\backslash}

\newcommand\SL{\mathrm{SL}}
\newcommand\PSL{\mathrm{PSL}}
\newcommand\SO{\mathrm{SO}}
\newcommand\PtZ{\PSL_2(\Z)}
\newcommand\PtR{\PSL_2(\R)}

\newcommand\Od{\mathcal{O}_d}
\newcommand\fraka{\mathfrak{a}}
\newcommand\ideal{\lhd}

\newtheorem{thm}{Theorem}[section]

\theoremstyle{definition}

\theoremstyle{remark}
\newtheorem*{rem}{Remark}

\newtheorem{example}{Example}[section]

\begin{document}

\title[Volumes of modular links]{Volumes of hyperbolic three-manifolds associated to modular links}

\author{Alex Brandts}
\address{McGill University, Montreal QC Canada}
\email{alexandre.brandts-longtin@mail.mcgill.ca}
\thanks{A.B.\ was supported by an NSERC USRA Fellowship}

\author{Tali Pinsky}
\address{Tata Institute for Fundamental Research, Mumbai India}
\email{tali@math.tifr.res.in}
\urladdr{http://www.math.tifr.res.in/~tali/}

\author{Lior Silberman}
\address{The University of British Columbia, Vancouver  BC  Canada}
\email{lior@math.ubc.ca}
\urladdr{https://www.math.ubc.ca/~lior/}
\thanks{L.S.\ was supported by an NSERC Discovery Grant}

\date{\today}

\subjclass[2010]{57M27; 11R29}

\begin{abstract}
Periodic geodesics on the modular surface correspond to periodic orbits
of the geodesic flow in its unit tangent bundle $\PtZ\bs\PtR$.
The complement of any finite number of orbits is a hyperbolic
$3$-manifold, which thus has a well-defined volume.
We present strong numerical evidence that, in the case of the set of geodesics
corresponding to the ideal class group of a real quadratic field, the volume
has linear asymptotics in terms of the total length of the geodesics.
This is not the case for general sets of geodesics.
\end{abstract}

\maketitle

\section{Introduction}
Let $X$ be a three-dimensional manifold.  A \emph{knot} in a $X$ is a simple
closed smooth curve, considered up to deformation; a \emph{link} in $X$ is
a collection of pairwise disjoint such curves, again considered up to
deformation.

We study numerically the topological type of knots and links
arising from closed geodesics on the modular surface;  These geodesics
have number-theoretic meaning, and we will show evidence for an
interesting behaviour for the links that arise when the geodesics are
grouped in a number-theoretic relevant way.

In the rest of this introduction we will review these general ideas,
discuss the motivation for our work, and outline the rest of the paper.

\subsection{Knots associated to closed geodesics}
Let $M$ be a Riemannian surface, that is a surface equipped with a Riemannian
metric.  Its \emph{unit tangent bundle} is the bundle whose fibre over
every $p\in M$ is the circle of direction vectors at $p$.  We note
that this bundle (denoted $T^1M$) is a three-dimensional manifold.

Now let $\gamma\colon [a,b]\to M$ be a smooth closed curve on $M$,
parametrized to have unit speed (smoothness at the endpoints requires
not merely that $\gamma(a) = \gamma(b)$ but that the derivatives match to
all orders).  Then for every $t\in [a,b]$ the derivative $\dot{\gamma}(t)$
is a unit (=direction) vector at $\gamma(t)\in M$, so that the map
$t\mapsto\left(\gamma(t),\dot{\gamma}(t)\right)$ describes a closed
curve in $T^1M$.  We require the self-intersection
of the curve to be \emph{transverse} (the direction vectors must be distinct)
so that the lift to $T^1M$ is not self-intersecting.  Included in this is
the assumption that the curve is \emph{primitive} (not a power of another
curve).
For such a curve, the (smooth) isotopy class of its lift in $T^1M$ is a knot therein.

Moreover, any smooth deformation of $\gamma$ into a curve $\gamma'$
can be modified to preserve the unit speed parametrization, at which point
it can be lifted to a deformation of the lifted curves in $T^1M$.  In other
words, the resulting knots depend only on the deformation (isotopy)
class of the curve.  In addition, we may construct a link by beginning with
several curves as long as all intersections (including self-intersections)
are \emph{transverse}.

Changing the metric on $M$ does not affect the topological properties of
$T^1M$ or the embeddings on curves, so (excluding the sphere and the torus,
which we do from now on) we may equip our surface with a finite-volume
metric of constant negative curvature $-1$, that is we realize it as a
hyperbolic surface (for the geometry of such surfaces see
\cite{Katok:FuchsianGroups}).
Fixing such a metric, each isotopy class of primitive
closed curves is of one of two types.  Generically, such a class
contains a unique representative of shortest length, necessarily a
closed geodesic.  If the surface is non-compact and hence has cusps
(infinitely long asymptotically narrow ends), there are in addition finitely
many classes of "periodic horocyles", curves going around a cusp.
Curves in the latter class can be deformed to be arbitrarily short by
``sliding them off the cusp" toward infinity.

\subsection{Modular links and their volumes}
One specific non-compact hyperbolic surface is the \emph{modular surface}
(a review of this surface appears in Section~\ref{sec:links}, we only note now that
the surface, its unit tangent bundle and its closed geodesics have explicit
realizations).

The recent surge of interest in "modular links", that is those that
arise from lifting closed geodesics on the modular surface, stems from
the work of Ghys (\cite{Ghys:KnotsDynamics}) where the closed geodesics
were connected to the periodic orbits of the Lorentz attractor, a
classical example of a chaotic dynamical system. 

Furthermore, there is a natural embedding of the unit tangent bundle
of the modular surface as the complement of a specific knot (the trefoil)
in the three-dimensional sphere.  Via this embedding we augment each modular
knot or link to a link in the three-sphere with an additional component
(this trefoil knot).  Ghys showed that the linking number of each modular knot
with the trefoil can be computed via the \emph{Rademacher function}, a
classical function on the modular group $\SL_2(\Z)$, and this allowed
Sarnak \cite{sarnak2010linking} to count geodesics by linking number.

Our contribution here is to study the topology of the augmented links
directly, specifically through properties of their complements.  Note that
the complement of the augmented link in the three-sphere is exactly
the complement of the original link in the unit tangent bundle.  For the
purposes of exposition it is more convenient to consider the complements
in the form $X\setminus\gamma$ where $X$ is the unit tangent bundle of the
modular surface and $\gamma$ is the (lift) of a family of closed
geodesics, while for computational purposes it will be convenient to
include the trefoil component in the calculations.

The starting point of our study is the fact that our link complements
$X\setminus\gamma$ always have a hyperbolic structure -- they can be
endowed with a finite-volume complete Riemann metric of constant
curvature $-1$ (for a proof see \cite{FoulonHasselblatt:ContactAnosov}).
By Mostow's Strong Rigidity Theorem this hyperbolic structure is
unique, so that metric properties of the resulting Riemannian $3$-manifold
are in fact topological properties of the link, depending only on the
deformation class of the curves making up $\gamma$.

We focus here on the simplest invariant of a hyperbolic three-manifold,
its volume, and specifically study the relationship between the volume
of $X\setminus \gamma$ and the complexity of $\gamma$ as measured by the
total length of its constituent geodesics (from the topological point
of view we essentially defined the length of a curve by taking the length
of the shortest representative in its deformation class).  Following
the usual convention in topology we shorten "the volume of the complement
of $\gamma$" to "the volume of $\gamma$".


In the related work sharing several authors
\cite{BergeronPinskySilberman:Volumes}, we have shown that there is
a universal constant $C$ such that the volume of a modular link is
bounded by $C$ times its length.  While a matching lower bound is
not possible in general (see Examples~\ref{exa:bounded-volume},~\ref{exa:log-volume} below), it is nevertheless
interesting to ask if our bound is sharp, that is if there are some links
whose volumes grow linearly with their length.  Experimentally this seems
to be the case for the links considered here.

\subsection{Results and structure of this paper}
Specifically, for each square-free positive integer $d>1$ let $\gamma_d$
be the set of geodesics corresponding to the ideal classes in the ring
of integers of the quadratic field $\QQ(\sqrt{d})$ (for details see
Section~\ref{sec:links}).  It is a beautiful theorem of Duke
\cite{Duke:SubConvexityEquidistribution} that these
sets of geodesics are becoming uniformly distributed on $X$: for any
compactly supported continuous function $f$ on $X$, the averages of $f$ on
the geodesics in $\gamma_d$ (with respect to the length measure)
coverge to the average of $f$ on $X$ (with respect to its Riemannian volume).
Inspired by this we examine numerically the volume of the complements
of these sets.

In Section \ref{sec:links} we discuss the modular surface, its unit
tangent bundle, and the
relationship between lifts of closed geodesics there and equivalence
classes of integral binary quadratic forms and construct the so-called
``fundamental" geodesics, those associated to the ideal classes
in the ring of integers of a real quadratic field, equivalently to 
classes of quadratic forms of fundamental discriminant.

In Section 3 we review the passage from number theory to topology:
the coding of geodesics by positive words in the two-letter alphabet
$\{X,Y\}$, how to obtain the coding of the a geodesic from the corresponding
quadratic form, as well as Williams's algorithm giving a combinatorial
description of the resulting link in $X$ from the list of words.

We use SageMath \cite{SageMath} and SnapPy \cite{SnapPy} to numerically
obtain the volume of the complement from this combinatorial description;
we give the details in Section 4.

Finally in Section 5 we present our main result (Figure 2):
for the family of links corresponding to the class groups of the
fields $\QQ\left(\sqrt{d}\right)$ with $d$ ranging over the square-free numbers
up to $1000$ (608 fields), the volumes seem to lie very close to the
best-fit line
$$ \textrm{volume} \approx 3.24 * \textrm{length} + 0.01\,.$$

\section{Modular Links and quadratic fields}\label{sec:links}

\subsection{Periodic geodesics on the modular surface}

Let $G = \PtR = \SL_2(\RR) / \left\{ \pm I\right\}$ be the quotient of 
the unimodular group $\SL_2(\RR)$ by its centre.  $G$ acts on the upper
half-plane $\HH = \{ x+iy \mid y>0  \}$ via the \emph{fractional linear
transformations}
$$\begin{pmatrix} a & b \\ c & d\end{pmatrix}\cdot z = \frac{az+b}{cz+d}\,.$$
The action is transitive and the stabilizer of $i$ is the circle group
$K = \SO(2)$, so that $\HH = G/K$.

The $G$ action on $\HH$ preserves the Riemannian metric
$ds^2 = \frac{dx^2+dy^2}{y^2}$, of constant negative curvature $-1$,
giving a realization of the hyperbolic plane.  The circle bundle $G\to \HH$
is exactly the unit tangent bundle.  Concretely,
the Gram--Schmidt decomposition of $g\in G$ (or Iwasawa decomposition) takes the form
$$g =	\begin{pmatrix} 1 & x\\ & 1\end{pmatrix}
	\begin{pmatrix} y^{1/2} & \\ & y^{-1/2}\end{pmatrix}
	\begin{pmatrix} \cos(\theta/2) & \sin(\theta/2) \\ -\sin(\theta/2) & \cos(\theta/2) \end{pmatrix},
$$
which we interpret as the tangent vector in direction $\theta$ at the point 
$z = x+iy \in \HH$ ($\theta=0$ corresponds to vectors pointing directly
upward in the half-plane realization above).

The \emph{geodesic flow} on the tangent bundle is the flow where the vector
$(z,\theta)$ flows along the geodesic starting at $z$ in the direction
$\theta$.  In terms of our representation above moving distance $t$ along 
the geodesic amounts to multiplying $g = (z,\theta)$ on the right by
$a_t = \begin{pmatrix} e^{t/2} & 0 \\ 0 & e^{-t/2}\end{pmatrix}$.

Now let $\Gamma = \PtZ$, a discrete subgroup of $G$.  The quotient
$M = \Gamma \bs \HH$ is the celebrated \emph{modular surface}.  Its points
parametrize, for example, lattices in $\R^2$ up to rotation and scaling.
Furthermore, the quotient space $X = \Gamma\bs G$ is clearly the
unit tangent bundle of $M$.  Geodesic flow on
$X = T^1M$ is still given by the right action of $A = \{ a_t \}_{t\in \R}$.

\begin{rem}
Since $\Gamma$ has elements of finite orders, the quotient $M$ is formally
an orbifold rather than a manifold.  This distinction, however, is
irrelevant to the discussions here.  These elements of finite order
will play a role in Section \ref{sec:coding}.
\end{rem}

In the realization $X$ of the tangent bundle, the geodesic through the
coset $x = \Gamma g\in X$ is periodic if and only if we have some
non-identity $a_t \in A$ such that
$x\cdot a_t = x$, equivalently such that $g\cdot a_t = \gamma g$.
Rewriting this as $g^{-1} \gamma g = a_t$ we see that a closed geodesic is
determined by a matrix $\gamma \in \Gamma$ having distinct real eigenvalues.
It is easy to check that conjugating $\gamma$ by $\eta\in\Gamma$ amounts
to replacing $g$ with $\eta g$, so that periodic geodesics are in bijection
with conjugacy classes in $\Gamma$ of elements diagonable over $\R$.
We also need to restrict to $\gamma$ which are not proper powers of other
elements of $\Gamma$ (``primitive" conjugacy classes) since powers amount to
geodesics going around the same loop multiple times.

For simplicity we shall use $\gamma$ to interchangeably denote the closed
geodesic on $M$, its lift to $X$, the corresponding primitive conjugacy class in
$\Gamma$ or a representative element of the conjugacy class.

\subsection{Quadratic forms and geodesics}
Each geodesic on $M$ is the projection to $M$ of a geodesic in $\HH$, and
such geodesics can be infinitely extended to have two endpoints on the
circle at infinity (this is similar to the case of straight lines in the flat
plane).  In the upper half-plane model the circle at infinity consists of the
real line (the obvious boundary of the upper half-plane) and a single point 
thought of as lying at imaginary infinity.  The geodesics accordingly take one
of two forms: either both of the endpoints lie on the real axis, and
then the geodesic is the Euclidean semicircle whose diameter is the interval
between the two endpoints, or one endpoint is at infinity, and then the
geodesic is the vertical line touching the other endpoint.

In particular a geodesic not passing through $i\infty$ can be encoded by
the quadratic polynomial whose roots are the two endpoints, equivalently
by any quadratic form $q(x,y)=ax^2 + bxy + cy^2$ such that the endpoints
are the roots of the polynomial $az^2 + bz +c$ (this is unique up to scaling)

It is well-known that the projection to $M$ of a geodesic is periodic exactly
when its endpoints are conjugate quadratic irrationalities, or equivalenty
when the quadratic form can be taken to have integral coefficients,
and we now study those.

An \emph{integral binary quadratic form} is the homogenous polynomial
$q(x,y) = ax^2 +bxy + cy^2$ with $a,b,c \in \Z$.  Its \emph{discriminant}
is the quantity $d = b^2-4ac$, and we will only consider the \emph{indefinite}
case where $d>0$.  Accordingly we will shorten ``indefinite integral binary
quadratic form" to ``quadratic form".

We give below a formula of a matrix $\gamma \in \Gamma$ representing the
conjugacy class corresponding to the geodesic determined by a quadratic form $q$.

Now $\Gamma$ acts on the quadratic forms by composition: if we think of $q$
as a function on $\Z^2$ then $q\circ\gamma$ is another quadratic form.
The orbits are called equivalence classes or simply \emph{classes} of
quadratic forms.
It is easy to verify that replacing $q$ with $q\circ\eta$
replaces the points $z_1,z_2$ with $\eta\cdot z_1,\eta\cdot z_2$ and hence
translates the geodesic $\gamma$ by $\eta$.  On the quotient we obtain the same
geodesic, and the result is a bijection between classes of quadratic forms
and periodic geodesics.

It is easy to check that equivalent forms have the same discriminant.
It was shown by Gauss that for each $d\neq 0$ there are finitely many
equivalence classes having discriminant $d$ (we concentrate on $d>0$).

It is useful to consider yet another number-theoretic description: each class
of quadratic forms with discriminant $d$ corresponds to a class of invertible
ideals in the quadratic order of discriminant $d$ (see for example the book
\cite{Cohn:AdvancedNumberTheory} which develops the theory of quadratic fields
through this correspondence, which is the historical point of view).
In the particular case of \emph{fundamental} discriminants (those which are
discriminants of quadratic fields), the correspondence is with ideals
classes in the ring of integers of the corresponding field and we give
one direction of the correspondence: the subring
$\Od = Z\left[\frac{d+\sqrt{d}}{2}\right]$ of the field
$K = \QQ\left(\sqrt{d}\right)$ is its \emph{ring of integers}.  Every ideal
$\fraka \ideal \Od$ is a free $\Z$-module of rank $2$.  Restricting the
\emph{norm form} $N^K_\Q$ to $\fraka$ therefore gives an integral quadratic
form.  Furthermore, equivalent ideals (those that differ by an element of
$K^\times$ correspond to equivalent quadratic forms, giving a bijection
between the group of ideal classes and the set of classes of quadratic
forms.

The group $\Od^\times$ of units (invertible elements) of the ring $\Od$
always has the form $\left\{\pm \epsilon^j\right\}_{j\in \ZZ}$, and we
call $\epsilon$ a \emph{fundamental unit} (there are four
such elements, but choosing any of them will give the same results).

Letting $m$ denote the squarefree part of $d$, we have
$K = \Q\left(\sqrt{m}\right)$, and we represent each element of $K$
in the form $x + \sqrt{m} y$.  Then $\Od$ contains the elements of this form
with $x,y\in \Z$ but may also contain those elements where $x,y$ are both
half-integers (elements in $\Z+\frac{1}{2}$), the latter case occurring exactly
when $m\equiv 1(4)$.  In this latter case we have $d=m$, and otherwise $d=4m$.

\section{Coding of geodesics}\label{sec:coding}

\subsection{Conjugacy classes in $\PtZ$ and words in two generators}\label{subsec:Conjugacy}
Continuing with our previous conventions we represent elements of $\PtZ$
by preimages in $\SL_2(\Z)$.  Specifically, 
let $U = \begin{pmatrix} 0 & 1\\ -1 & 0\end{pmatrix}$
and $V = \begin{pmatrix} 0 & -1\\ 1 & 1\end{pmatrix}$.
It is not hard to check that $U,V\in\PtZ$ are elements of order $2,3$
respectively, and it is well-known that $\PtZ$ is the free product of the
cyclic groups generated by $U$ and $V$.  Accordingly every element of the group
has a unique representation as a \emph{reduced} word in $U,V$, that is a
word in this alphabet containing no subword of the form $UU$ or $VVV$
(note that only positive words are necessary since $U^{-1} = U$ and
$V^{-1}=V^2$.

Recalling that periodic geodesics correspond to (certain)
\emph{conjugacy classes} in $\PtZ$, we now discuss canonical representatives
for those.  Note that if a word beings and ends in $U$ then it can be thought
of as beginning with $U$ and ending with $U^{-1} = U$ so deleting the $U$s
at the ends doesn't affect the conjugacy class.  Similarly if the words beings
with $V$ or $V^2$ one can conjugate by $V^2$ or $V$ to obtain a word
beginning with $U$.  The elements of finite order in $\PtZ$ are exactly those
conjugate to $1,U,V,V^2$, and we see that every conjugacy class of elements
of infinite order contains a shortest representative which is a positive word
in the alphabet $X,Y$ where $X=UV$, $Y=UV^2$.  Note that
$X = \begin{pmatrix} 1 & 1\\ 0 & 1\end{pmatrix}$
and $Y = \begin{pmatrix} 1 & 0\\ 1 & 1\end{pmatrix}$.

The shortest representative of the conjugacy class is unique (as a word in $U,V$)
up to cyclic permutation (this is a general fact about free products).  It follows
that every conjugacy class of elements of infinite order has a representative
as a word in $X,Y$, unique up to cyclic permutation.  We also know that an element
of $\PtZ$ is a power of another if and only if the corresponding word in $X,Y$
is periodic, so that primitive conjugacy classes correspond to primitive words
(those that don't consist of repetitions of a subword).  We also note that every
parabolic element of $\PtZ$ is conjugate to a power of $X$ (possibly negative).
Since $Y$ is conjugate to $X^{-1}$ we see that the words of the form $X^k$
or $Y^k$ ($k\geq 1$) exactly parametrize the parabolic conjugacy classes.

Conclusion: primitive closed geodesics on the modular surface are in bijection
with positive words in $X,Y$ containing both symbols.

\subsection{From matrices to words}\label{sec:reduction}
For sake of completeness, we mention the algorithmic solution to
obtaining the word in the generators $X,Y$, given the matrix representation of a hyperbolic
element of $\PtZ$.

Given a matrix $A=\begin{pmatrix} a & b\\ c & d\end{pmatrix}\in\PtZ$, we first represent $A$ as a product of $U$ and $X=UV$, as follows:
Multiplying any matrix $A$ by $U$ will switch its first and second rows, and multiply the second row by $-1$.
Next, multiplying by $(UV)^n$ will add $n$ times the second row to the first.
Thus, we may assume that $|a|\geq|c|$ by multiplying by $U$ if needed. We may also assume that $a$ is positive as $A\equiv-A$.
Then we can subtract $n$ times $|c|$ from $a$, for the largest $n$ such that $a-nc\geq0$. Continuing in this way, we are actually applying the Euclidean algorithm to $a$ and $c$. As $\det(A)=ad-bc=1$, $a$ and $c$ are relatively prime, and thus in a finite number of steps we reach a matrix
$A'=\begin{pmatrix} 1 & 0\\ c & d\end{pmatrix}$. 
But as $A'\in\PtZ$ $d=1$, and so $A'=(UV)^m$ for some integer $m\neq0$.
This results in a decomposition of $A$ into the generators $U$ and $V$. By switching to positive powers and grouping the generators into subsequences of $UV$ and $UV^2$, we reach the desired representation.

\subsection{Curves from words and the template}
In order to investigate the topology of the closed geodesic we need a way
to explicitly construct a representative of the isotopy class of the geodesic
on the unit tangent bundle from a  representation of the conjugacy class.
This was done by Ghys \cite{Ghys:KnotsDynamics}, as follows.

A \emph{template} \cite{BirmanWilliams2}  is a branched surface with boundary and an expansive semi flow.
The following is a theorem of Ghys, relying on a theorem of Birman and Williams \cite{BirmanWilliams2} which asserts that any Anosov flow on a three manifold has a template.

\begin{figure}[ht]
\begin{overpic}[width=6cm]{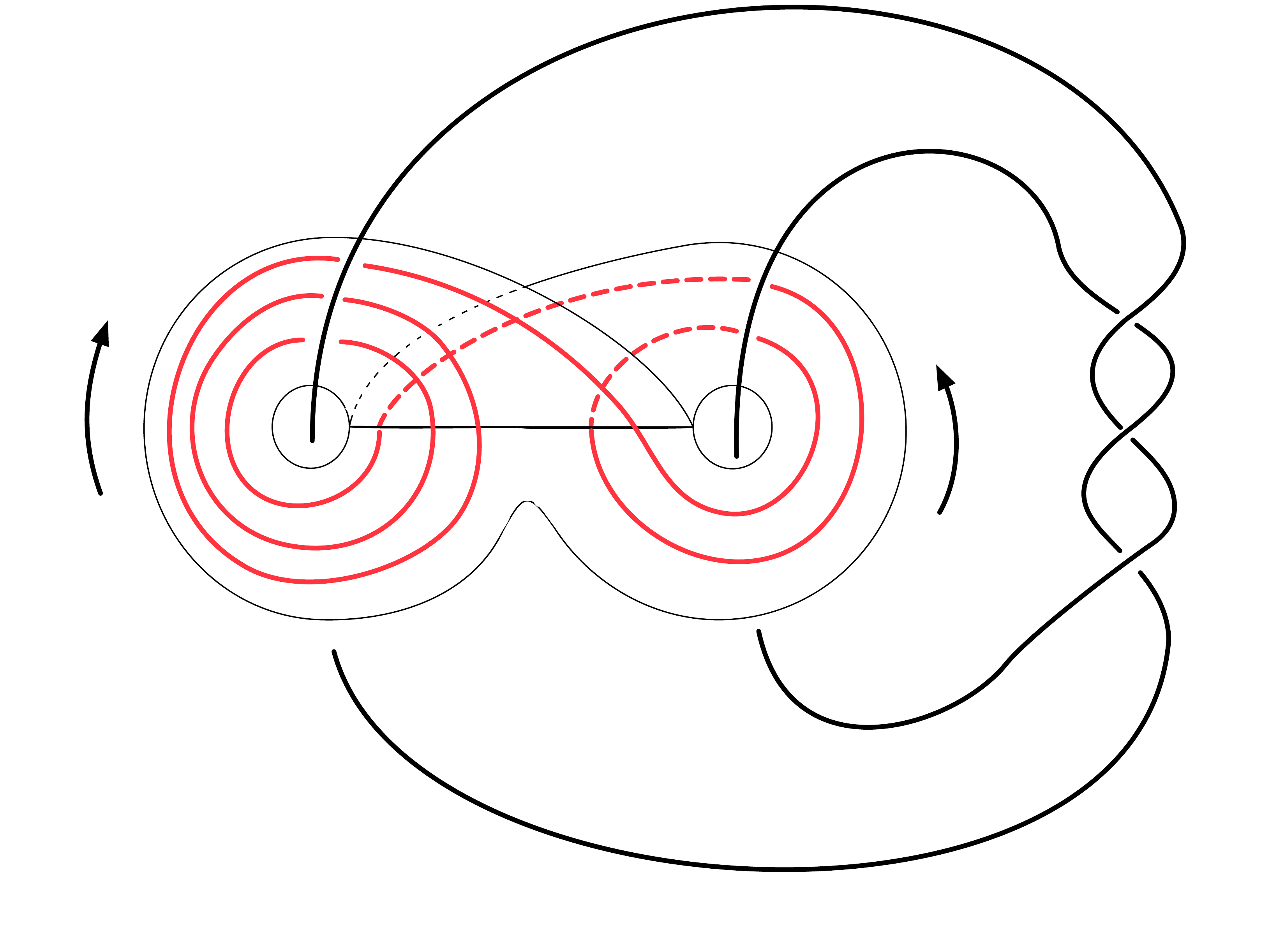}
\put(1,40){$x$}
\put(77,40){$y$}
\end{overpic}
\caption{The template $T$ for the geodesic flow on $\MS$, together with a periodic orbit. The arrows indicate the direction of the semiflow.}\label{fig:template}
\end{figure}

\begin{thm} [\cite {Ghys:KnotsDynamics}] The set of periodic geodesics on $\MS$ is in bijective correspondence with the set of periodic orbits on the template $T$ embedded in $S^3\setminus\text{trefoil}$ as depicted in figure \ref{fig:template}, excluding the boundary curves of $T$.
On any finite subset, the correspondence is by an ambient isotopy.
\end{thm}

This means that in order to study knot properties of closed modular geodesics, one may study closed orbits on the template $T$.
The template comes with a symbolic dynamics given by the symbols $x$ and $y$ that correspond to passing through the left or through the right ear (or equivalently through the left or the right half of the branch line, which is a Birkhoff section for the semiflow).
Ghys proves in \cite{Ghys:KnotsDynamics} that the symbolic dynamics of the representation  of a geodesic in $T$ is equivalent to its representation as a word in the generators $X,Y$ for $\PtZ$.


\section{Computational method}

We start with some table of quadratic fields and generate the matrices as explained in Section \ref{subsec:Conjugacy}.

To decompose a matrix into a product of $X$ and $Y$, we use a standard algorithm to decompose matrices in $SL_2(\Z)$.

Then, given the corresponding periodic words in $x,y$, Williams \cite{Williams} provides an algorithm for locating the periodic orbit corresponding to it on $T$, as follows.
 Start by listing the word $w$ and all its cyclic permutations as $w_1=w,\cdots,w_k$. Each cyclic permutation $w_j$ corresponds to a point $p_j$ in the intersection of the orbit with the branchline, where the word $w_j$ is the orbit when you start to read it from $p_j$.
 Consider the lexicographic order induced by $x<y$. The position of the points $p_j$ along the branchline is given by the lexicographic order. i.e., if $p_j<p_k$ then $p_j$ is further to the left in Figure ~\ref{fig:template}. Thus one orders the points $p_1,\dots,p_k$, and connects each $p_j$ to $p_{j+1}$ in the obvious way given that flowlines do not intersect.

This yields a way to convert the words into a so called DT-code for the link in $S^3$ given by the union of the geodesics in $\gamma$ and the trefoil knot:
Starting at the point $p_1$, we trace the geodesic. It starts through the left ear, therefore we have the first crossing is an overcrossing  (and we assign a 1) with the trefoil knot and then the second is an undercrossing (assigned a $-2$). Recall we keep track of the sign of a crossing only for an even number. All the crossings above the branch line have the same orientation, i.e., coming from the $x$ ear they will all be overcrossings. 
Suppose we reach a point $p_j$ (recall we keep track of what point connects to what point when generating the order of these points). Then we will have $j-1$ crossings above the branchline (assigned the positive numbers 3 to $j+2$). We continue in this way, passing to the next leftmost component of $\tilde\gamma$ once we reached $p_1$ and finally tracing the trefoil knot once we have traced all components in $\tilde\gamma$.

We then feed the resulting DT code into SnapPy to obtain an estimate of the volume.

To compute the sum of the geometric lengths, let $\epsilon$ be the fundamental unit, $R = |\log\epsilon|$ the regulator, $h$ the class number, $h^+$ the narrow class number.

There are two cases:

1. The fundamental unit $\epsilon$ has norm $+1$.  In this case $|\log\epsilon|$
  is the regulator and the narrow class number is twice the class
  number.  Thus the sum of the geometric lengths is given by $h^+\cdot |\log\epsilon| = 2hR$.

2. The fundamental unit $\epsilon$ has norm $-1$.  In this case the unit we use is
  not the fundamental unit, but $\epsilon^2$. i.e., its log is twice the regulator. On the other hand, the narrow class number is  in this case the same as the class number, and thus the sum of the geometric lengths is given by $h^+ \cdot |\log\epsilon^2| = 2hR$ in this case as well.

\section{Results and discussion}
\begin{figure}[ht]
\begin{overpic}[width=15cm]{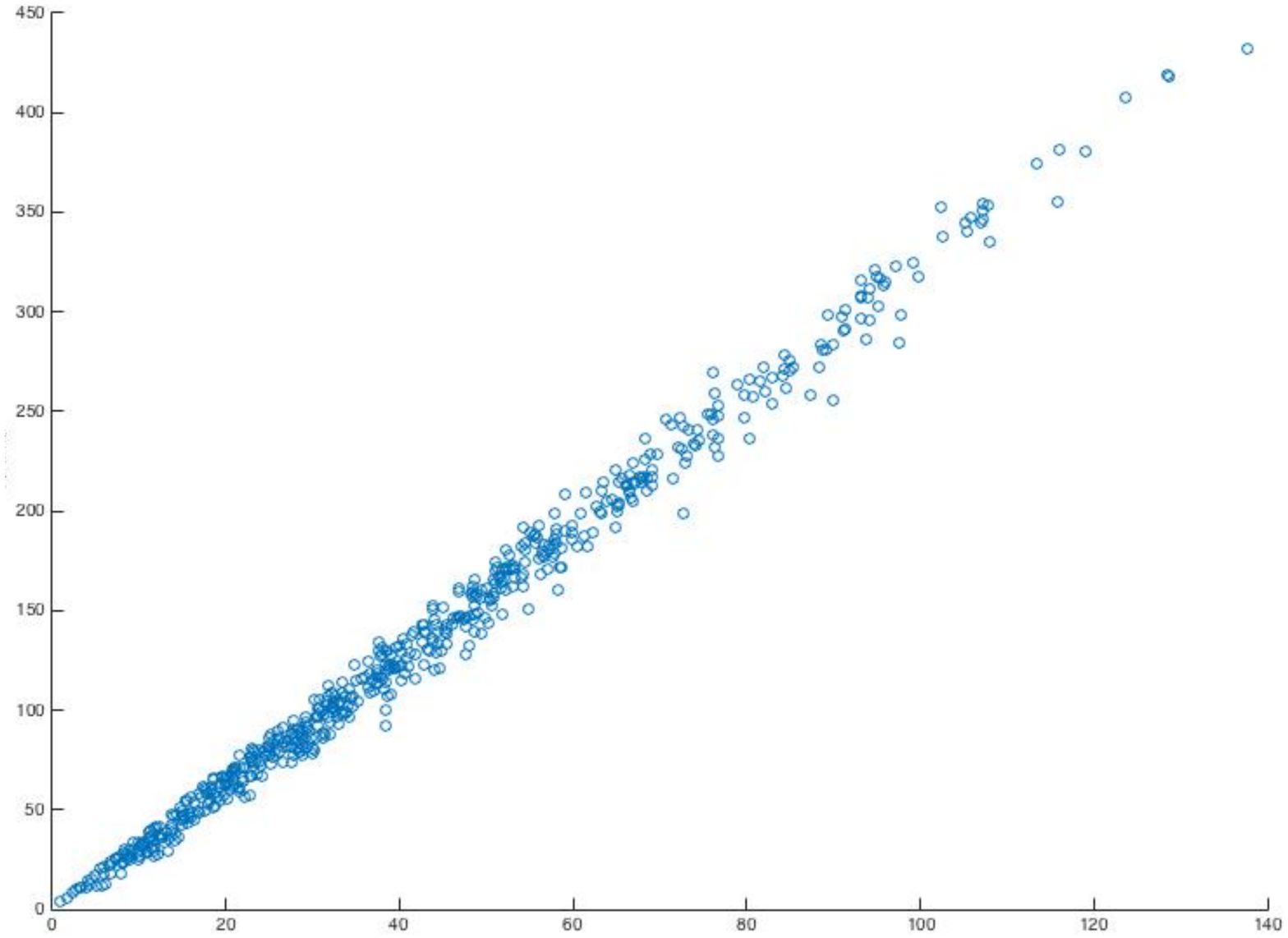}
\put(40,3){geometric length}
\put(5,35){ \begin{sideways} volume \end{sideways}}
\end{overpic}
\caption{}\label{fig:growth}
\end{figure}

The volumes of the complements versus the geometric length computed as in the previous section are shown for discriminants up to 3992 (or to square free integers up to 1000) in Figure ~\ref{fig:growth}
and we find that the growth is linear. 

In contrast, there are families of geodesics for which the lengths tend to infinity while the volume is bounded by a constant:

\begin{example}\label{exa:bounded-volume}
Let $\gamma$ be a geodesic on $\MS$ coded by the word $x^ny^m$ for some natural numbers $n$ and $m$.
Its length is roughly proportional to $\log(n+m)$ and, in particular, tends
to infinity with $n$ and $m$.  The number of crossings on the template (which will be proportional  to the number of tetrahedra in our triangulation) is then $n+m$.
On the other hand, these geodesics correspond to a knot winding more and more
around the trefoil, first in one template ear and then in the other.
Thus, their volumes are all bounded by the volume of
$T^1\MS\setminus(\gamma\cup\alpha\cup\beta)$ where $\gamma$ is the geodesic
corresponding to $xy$, and $\alpha$ and $\beta$ are both trivial knots
encircling one strand of $\gamma$ and the strand of the trefoil in the
centre of the corresponding ear (c.f Adams \cite{adams1985thrice}). 
\end{example}

Note that the period is equal to 1 for any geodesic $\gamma$ in the above example.
Initially, the authors believed that the period would provide a lower bound for the volume, and thus a family shown to efficiently increase the period of the continued fraction expansion would have linear growth. However, this is not the case:

\begin{example}\label{exa:log-volume}
Consider the family of geodesics with $w(\gamma_n)=x(xy)^n$ on $\MS$. The length of $\gamma_n$ grows linearly in the period $n$, however the volumes are bounded. The fact that it can grow at most logarithmically follows from the arguments given in \cite{BergeronPinskySilberman:Volumes}, and numerically we find that the sequence of volumes has a limit $\approx15$.
\end{example}

It is interesting to note that the authors have not been able to find any family of geodesics that obtains a linear growth without relying on the class field.
\bibliographystyle{alpha}
\bibliography{numerics.bib}

\end{document}